\documentclass[12pt,a4paper,twoside]{amsart} 
\usepackage{amssymb,amsbsy,amsmath,amsfonts,amssymb,amscd,times,
graphics,color,xypic,footmisc,fancyhdr,multicol,fancybox,
graphicx,mathrsfs,rotating,ifthen,wasysym}
\usepackage{times,color} 

\newcommand{\HEAD}[2]{%
\pagestyle{fancy}
\fancyhead[RO]{\scriptsize\sf\thepage}
\fancyhead[CO]{{\scriptsize\sf #1}}
\fancyhead[LE]{\scriptsize\sf\thepage}
\fancyhead[CE]{{\scriptsize\sf #2}}
\fancyfoot{}}

\let\mathcal\mathscr 

\newtheorem{The}{Theorem}[section] 
\newtheorem{Theorem}{Theorem}[section] 
\newtheorem{Proposition}[The]{Proposition}

\newtheorem{Observation}[The]{Observation} 

\theoremstyle{definition}

\newcommand{\C}{\mathbb{C}}\newcommand{\N}{\mathbb{N}}

\newcommand{\R}{\mathbb{R}} 

\newcommand{\vf}{\vfill\end{document}}

\begin{document} 

%\large

$\:$

\bigskip\bigskip

\begin{center}

{\large\bf On Transfer of Biholomorphisms}

\smallskip

{\large\bf Across Nonminimal Loci}

\end{center}

\bigskip

\begin{center}
Jo\"el {\sc Merker}
\end{center}

\medskip

\begin{center}
\begin{minipage}[t]{10.25cm}
\baselineskip =0.32cm 
{\scriptsize
{\bf Abstract.}
A connected real analytic hypersurface $M \subset \C^{ n+1}$
whose Levi form is nondegenerate in at least one point\,\,---\,\,hence 
at every point of some Zariski-dense open subset\,\,---\,\,is 
locally biholomorphic to the model Heisenberg quadric pseudosphere
of signature $(k, n - k)$ in one point 
if and only if, at every other Levi nondegenerate point, it
is also locally biholomorphic to some Heisenberg
pseudosphere, possibly having a different signature $(l, n-l)$.
Up to signature, pseudo-sphericity then jumps across the 
Levi degenerate locus, and in particular, across the nonminimal
locus, if there exists any.
}
\end{minipage}
\end{center}

\bigskip

\begin{center}
\begin{minipage}[t]{11.75cm}
\baselineskip =0.35cm {\scriptsize

\centerline{\bf Table of contents}

\smallskip

{\bf \ref{introduction}.~Introduction
\dotfill~\pageref{introduction}.}

{\bf \ref{proof-C-2}.~Proof in $\C^2$
\dotfill~\pageref{proof-C-2}.}

{\bf \ref{C-n-1}.~Proof in $\C^{n+1}$ ($n \geqslant 2$)
\dotfill~\pageref{C-n-1}.}

}\end{minipage}
\end{center}

\bigskip

\medskip

\date{\number\year-\number\month-\number\day}

\begin{abstract} 
zero curvature equations
\end{abstract} 

\pagestyle{headings} 
\markright{Pseudo-spherical nonminimal hypersurfaces} 

\section{Introduction} 
\label{introduction}
\HEAD{1.~Introduction}{
Jo\"el Merker, D\'epartement de Math\'ematiques d'Orsay, Paris, France}

\noindent

The goal of this paper is to provide the complete details of an
alternative direct proof of a recent theorem due to Kossovskiy and
Shafikov (\cite{ Kossovskiy-Shafikov-2012})
which relies on the explicit zero-curvature equations
obtained in~\cite{ Merker-2010, Merker-2009}, following the lines of a
clever suggestion of Beloshapka. In fact, 
the proof we give here freely brings a more general statement.

Let $M \subset \C^{ n+1}$ be a connected real analytic hypersurface
with $n \geqslant 1$. One says that $M$ is $(k, n-k)$ {\sl
pseudo-spherical} at one of its points $p$ if it is locally near $p$
biholomorphic to some Heisenberg $(k, n-k)$-pseudo-sphere
having, in coordinates $(z_1, \dots, z_n, w) \in
\C^{ n+1}$, the model 
quadric equation:
\[
w
=
\overline{w}
+
2i\big(
-\,z_1\overline{z}_1
-\cdots-
z_k\overline{z}_k
+
z_{k+1}\overline{z}_{k+1}
+\cdots+
z_n\overline{z}_n
\big),
\] 
for some integer $k$ with $0 \leqslant k \leqslant n - k$;
when $n = 1$, one simply says that $M$ is {\sl spherical}.

It is known that a connected real analytic hypersurface of $\C^{ n+1}$
which is Levi nondegenerate at every point is
$(k,n-k)$-pseudospherical at one point if and only if it is $(k,
n-k)$-pseudospherical at every point. More generally, 
we establish that propagation of $(k, n-k)$
pseudo-sphericality also holds in presence of Levi degenerate points
of arbitrary kind.

\begin{Theorem}
\label{main-theorem}
Let $M \subset \C^{ n+1}$ be a connected real analytic 
geometrically smooth hypersurface which is Levi nondegenerate
in at least one point (hence in some nonempty open subset). Then:

\medskip\noindent{\bf (a)}
The set of Levi nondegenerate points of $M$ is a Zariski open subset
of $M$ in the sense that there exists a certain {\em
proper}\,\,---\,\,{\em i.e.} having dimension $\leqslant \dim\, M -
1 = 2n$\,\,---\,\,locally closed global 
real analytic subset $\Sigma_{\sf LD}
\subset M$ locating exactly the Levi degenerate points of $M$:
\[
p\,\in\,M\big\backslash\Sigma_{\sf LD}
\,\,\,\Longleftrightarrow\,\,\,
\text{\small\sf Levi form of $M$ at $p$ is nondegenerate}.
\]

\medskip\noindent{\bf (b)}
If $M$ is locally biholomorphic, in a neighborhood of one of its
points $p \in M$, to some Heisenberg $(k, n-k)$-pseudo-sphere having
equation:
\[
w
=
\overline{w}
+
2i\big(
-\,z_1\overline{z}_1
-\cdots-
z_k\overline{z}_k
+
z_{k+1}\overline{z}_{k+1}
+\cdots+
z_n\overline{z}_n
\big),
\] 
for some integer $k$ with $0 \leqslant k \leqslant n - k$ (so that $p
\in M \big\backslash \Sigma_{\sf LD}$ necessarily is a Levi
nondegenerate point of $M$ too), then locally at every other Levi
nondegenerate point $q \in M \big\backslash \Sigma_{\sf LD}$, 
the hypersurface $M$ is also locally biholomorphic to some
Heisenberg $(l, n - l)$-pseudo-sphere, with, possibly $l \neq k$.
\end{Theorem}

Surprisingly, Example 6.2 in~\cite{ Kossovskiy-Shafikov-2012} shows
that $l \neq k$ may occur, in the case of a {\sl nonminimal}
hypersurface of $\C^{ n+1}$ with $n \geqslant 2$, a local example for
which the Levi degenerate locus $\Sigma_{\sf LD}$ consists precisely
of a complex $n$-dimensional hypersurface contained in $M$.

\medskip\noindent{\bf Acknowledgments.}
The author is grateful to Valerii Beloshapka, to Ilya Kossovskiy, to
Rasul Shafikov, and also to an objective benevolent anonymous referee
in the field of CR mappings for raising the interest of making public
the details of an alternative proof of the theorem ({\em cf.} 
a few e-mail
exchanges in March-April 2013 during
which the easy proof presented here was explained). 
He is also grateful to Alexander
Sukhov for introducing him to the subject of the interactions between
CR geometry and partial differential equations ({\em cf.} \cite{
Sukhov-2001, Sukhov-2002, Sukhov-2003}). Lastly, the above theorem
was mentioned in a talk of the author at the {\sl Abel Symposium 2013},
organized by John-Erik {\sc Forn{\ae}ss}, Marius {\sc Irgens} and
Erlend {\sc Forn{\ae}ss-Wold} ({\em cf.}~\cite{ Merker-Abel-2013}), a
talk during which an expert in nonminimal hypersurfaces was absent.

\section{Proof in $\C^2$} 
\label{proof-C-2}
\HEAD{\ref{proof-C-2}.~Proof in $\C^2$}{
Jo\"el Merker, D\'epartement de Math\'ematiques d'Orsay, Paris, France}

Let $M \subset \C^2$ be a connected real analytic
hypersurface. Pick a point:
\[
p\in M,
\]
and choose some affine coordinates centered at $p$: 
\[
(z,w)
=
\big(x+iy,\,u+iv\big)
\]
satisfying:
\[
T_0
M
=
\big\{{\sf Im}\,w=0\big\},
\]
so that the implicit function theorem represents $M$ as:
\[
u
=
\varphi(x,y,v),
\]
in terms of some {\sl graphing function} $\varphi$ which
is expandable in convergent Taylor series in some (possibly small)
open bidisc:
\[
\square_{\rho_0}^2
\,:=\,
\big\{
(z,w)\in\C^2\colon\,
\vert z\vert<\rho_0,\,\,
\vert w\vert<\rho_0
\big\},
\]
with of course $\rho_0 > 0$.

Classically, writing:
\[
{\textstyle{\frac{w+\overline{w}}{2}}}
=
\varphi
\big(
{\textstyle{\frac{z+\overline{z}}{2}}},\,
{\textstyle{\frac{z-\overline{z}}{2i}}},\,
{\textstyle{\frac{w-\overline{w}}{2i}}}
\big)
\]
and using the analytic implicit function theorem, 
one solves $w$ in terms of $z, \overline{ z}, 
\overline{ w}$ getting a representation of
$M$ as:
\[
w
=
\Theta\big(z,\overline{z},\overline{w}\big);
\]
recall that implicitly, when one does this, one must considers
$(z, w, \overline{ z}, \overline{ w})$ as 
$4$ independent complex variables, which amounts
to {\sl complexify} them, namely to introduce
the {\sl complexified variables}:
\[
\big(z,w,\underline{z},\underline{w}\big)
\in
\C^4;
\]
in what follows, we will work with 
$(z, w, \overline{ z}, \overline{ w})$-variables, 
keeping in mind that they can be replaced by $(z, w,
\underline{ z}, \underline{ w})$ {\em since all objects are convergent
Taylor series}; so here, $\Theta ( z, \underline{ z}, \underline{ w}
)$ is a convergent Taylor series of $(z, \underline{ z}, \underline{
w})$ for $\vert z \vert < \rho_0$, $\vert \underline{ z} \vert <
\rho_0$, $\vert \underline{ w} \vert < \rho_0$, after shrinking
$\rho_0 > 0$ if necessary.

Moreover, since:
\[
0
=
\varphi(0)
=
\varphi_x(0),
=
\varphi_y(0)
=
\varphi_u(0),
\]
one has:
\[
\Theta
=
\overline{w}
+
{\rm O}(2).
\]

Now, it is known\,\,---\,\,or it could be taken here
as a definition\,\,---\,\,that $M$ is {\sl Levi nondegenerate}
at $0 \in M$ when the local holomorphic map:
\[
\aligned
\C^2
&
\,\longrightarrow\,\,
\C^2
\\
\big(\overline{z},\overline{w}\big)
&
\,\longmapsto\,
\Big(
\Theta\big(0,\overline{z},\overline{w}\big),\,\,
\Theta_z\big(0,\overline{z},\overline{w}\big)
\Big)
\endaligned
\]
is of rank $2$ at $( \overline{z}, \overline{w}) = (0, 0)$;
of course, one would better think in
terms of $( \underline{ z}, \underline{ w})$-variables
here.

More generally, $M$ is 
{\sl Levi nondegenerate} at an arbitrary point
close to the origin:
\[
(z_p,w_p)
\in
M
\cap 
\square_{\rho_0}^2
\]
when the map:
\[
\aligned
\C^2
&
\,\longrightarrow\,\,
\C^2
\\
\big(\overline{z},\overline{w}\big)
&
\,\longmapsto\,
\Big(
\Theta\big(z_p,\overline{z},\overline{w}\big),\,\,
\Theta_z\big(z_p,\overline{z},\overline{w}\big)
\Big)
\endaligned
\]
is of rank $2$ at $(\overline{ z}_p, 
\overline{ w}_p)$, which precisely means
the nonvanising of the Jacobian determinant:
\[
\aligned
0
\neq
\det
\left(\!
\begin{array}{cc}
\Theta_{\overline{z}} & \Theta_{\overline{w}}
\\
\Theta_{z\overline{z}} & \Theta_{z\overline{w}}
\end{array}
\!\right)
= 
\Theta_{\overline{z}}\,
\Theta_{z\overline{w}}
- 
\Theta_{\overline{w}}\,
\Theta_{z\overline{z}}
\endaligned
\]
at:
\[
(z,\overline{z},\overline{w})
=
(z_p,\overline{z}_p,\overline{w}_p).
\]

One may either show-check that such a definition 
re-gives the standard definition of 
Levi nondegeneracy ({\em cf.}~\cite{ 
Merker-2005a, Merker-2005b, Merker-Porten-2006}), 
or prove directly that as it stands, it
really is independent of the choice of coordinates
(\cite{ Merker-Pocchiola-Sabzevari-2013}).

Although we could then spend time to re-prove it properly, 
the following
fact\,\,---\,\,here admitted\,\,---\,\,is well known.

\begin{Proposition}
If a connected real analytic hypersurface $M^{ 2n+1} 
\subset \C^{ n + 1}$ is Levi
nondegenerate in at least one point, then the set of Levi degenerate
points of $M$ is a proper real analytic subset: 
\[
\Sigma_{\sf LD}
\subsetneqq 
M.
\qed
\]
\end{Proposition}

In Theorem~\ref{main-theorem}, we indeed make the assumption that 
$\Sigma_{ \sf LD}$ is proper, since otherwise, the real analytic
$M$ would be {\sl Levi-flat}, hence as is known, everywhere locally
biholomorphic to $\C^n \times \R$.

\medskip

Suppose to begin with for $M^3 \subset \C^2$ that:
\[
0
\not\in
\Sigma_{\sf LD}.
\]
Then the above map being of rank $2$ at
$(\overline{ z}, \overline{ w}) = (0,0)$, 
one can solve, following~\cite{
Merker-2010}, the following two equations:
\[
\aligned
w(z)
&
=
\Theta\big(z,\overline{z},\overline{w}\big),
\\
w_z(z)
&
=
\Theta_z\big(z,\overline{z},\overline{w}\big),
\endaligned
\]
for the two variables $(\overline{ z}, \overline{ w})$, and 
then insert the latter in:
\[
w_{zz}(z)
=
\Theta_{zz}\big(z,\overline{z},\overline{w}\big),
\]
to get a complex second-order ordinary differential equation:
\[
w_{zz}(z)
=
\Phi\big(z,w(z),w_z(z)\big).
\]
One should notice that the possibility
of solving $(\overline{ z}, \overline{ w})$
is expressed by the nonvanishing
of {\em exactly and precisely the same} Jacobian determinant as the
one expressing Levi nondegeneracy. 

In the article~\cite{ Merker-2010}, one deduces from an 
explicitly known condition on $\Phi$ for this second-order equation
$w_{ zz} = \Phi( z, w, w_z)$ 
to be pointwise equivalent to the free particle Newtonian equation:
\[
w_{z'z'}'
=
0,
\]
that a real analytic hypersurface $M \subset \C^2$ 
which is Levi nondegenerate at $0 \in M$ as above
is {\sl spherical} in the sense\,\,---\,\,recall the 
definition\,\,---\,\,that it is
locally biholomorphic to:
\[
w'
=
\overline{w}'
+
2i\,z'\overline{z}',
\]
if and only if its complex graphing function $\Theta$
satisfies an explicit (not completely developed) equation
which we now present.

Introduce the expression:

\[
\footnotesize
\aligned
{\sf AJ}^4(\Theta)
&
:=
\frac{1}{
[\Theta_{\overline{z}}\Theta_{z\overline{w}}
-\Theta_{\overline{w}}\Theta_{z\overline{z}}]^3}
\bigg\{
\Theta_{zz\overline{z}\overline{z}}
\bigg(
\Theta_{\overline{w}}\Theta_{\overline{w}}
\left\vert\!\!
\begin{array}{cc}
\Theta_{\overline{z}} & \Theta_{\overline{w}}
\\
\Theta_{z\overline{z}} & \Theta_{z\overline{w}}
\end{array}
\!\!\right\vert
\bigg)
-
\\
&
\ \ \ \ \
-\,
2\Theta_{zz\overline{z}\overline{w}}
\bigg(
\Theta_{\overline{z}}\Theta_{\overline{w}}
\left\vert\!\!
\begin{array}{cc}
\Theta_{\overline{z}} & \Theta_{\overline{w}}
\\
\Theta_{z\overline{z}} & \Theta_{z\overline{w}}
\end{array}
\!\!\right\vert
\bigg)
+
\Theta_{zz\overline{w}\overline{w}}
\bigg(
\Theta_{\overline{z}}\Theta_{\overline{z}}
\left\vert\!\!
\begin{array}{cc}
\Theta_{\overline{z}} & \Theta_{\overline{w}}
\\
\Theta_{z\overline{z}} & \Theta_{z\overline{w}}
\end{array}
\!\!\right\vert
\bigg)
+
\\
&
\ \ \ \ \
+
\Theta_{zz\overline{z}}
\bigg(
\Theta_{\overline{z}}\Theta_{\overline{z}}
\left\vert\!\!
\begin{array}{cc}
\Theta_{\overline{w}} & \Theta_{\overline{w}\overline{w}}
\\
\Theta_{z\overline{w}} & \Theta_{z\overline{w}\overline{w}}
\end{array}
\!\!\right\vert
-
2\Theta_{\overline{z}}\Theta_{\overline{w}}
\left\vert\!\!
\begin{array}{cc}
\Theta_{\overline{w}} & \Theta_{\overline{z}\overline{w}}
\\
\Theta_{z\overline{w}} & \Theta_{z\overline{z}\overline{w}}
\end{array}
\!\!\right\vert
+
\Theta_{\overline{w}}\Theta_{\overline{w}}
\left\vert\!\!
\begin{array}{cc}
\Theta_{\overline{w}} & \Theta_{\overline{z}\overline{z}}
\\
\Theta_{z\overline{w}} & \Theta_{z\overline{z}\overline{z}}
\end{array}
\!\!\right\vert
\bigg)
+\\
&
\ \ \ \ \
+
\Theta_{zz\overline{w}}
\bigg(
-\Theta_{\overline{z}}\Theta_{\overline{z}}
\left\vert\!\!
\begin{array}{cc}
\Theta_{\overline{z}} & \Theta_{\overline{w}\overline{w}}
\\
\Theta_{z\overline{z}} & \Theta_{z\overline{w}\overline{w}}
\end{array}
\!\!\right\vert
+
2\Theta_{\overline{z}}\Theta_{\overline{w}}
\left\vert\!\!
\begin{array}{cc}
\Theta_{\overline{z}} & \Theta_{\overline{z}\overline{w}}
\\
\Theta_{z\overline{z}} & \Theta_{z\overline{z}\overline{w}}
\end{array}
\!\!\right\vert
-
\Theta_{\overline{w}}\Theta_{\overline{w}}
\left\vert\!\!
\begin{array}{cc}
\Theta_{\overline{z}} & \Theta_{\overline{z}\overline{z}}
\\
\Theta_{z\overline{z}} & \Theta_{z\overline{z}\overline{z}}
\end{array}
\!\!\right\vert
\bigg)
\bigg\},
\endaligned
\]
noticing that its denominator: 
\[
\big[\Theta_{\overline{z}}\Theta_{z\overline{w}}
-\Theta_{\overline{w}}\Theta_{z\overline{z}}\big]^3
\]
does not vanish at the origin
since $0 \in M$ was assumed (temporarily) to be a 
Levi nondegenerate point.
Introduce also the vector field derivation:
\[
\mathcal{D}
:=
\frac{-\,\Theta_{\overline{w}}}{
\Theta_{\overline{z}}\Theta_{z\overline{w}}
-\Theta_{\overline{w}}\Theta_{z\overline{z}}}\,
\frac{\partial}{\partial\overline{z}}
+
\frac{\Theta_{\overline{z}}}{
\Theta_{\overline{z}}\Theta_{z\overline{w}}
-\Theta_{\overline{w}}\Theta_{z\overline{z}}}\,
\frac{\partial}{\partial\overline{w}}.
\]
Then the main and unique theorem of~\cite{ Merker-2010}
states that $M$ is spherical at $0$ if and only if:
\[
0
\equiv
\mathcal{D}\big(\mathcal{D}\big({\sf AJ}^4(\Theta)\big)\big),
\]
identically in $\C\big\{ z, \overline{ z}, \overline{ w} \big\}$.

Unfortunately, it is essentially impossible to print in 
a published article what one
obtains after a full expansion of these two derivations.

Nevertheless, by thinking a bit, one convinces oneself that
after full expansion, and reduction to a common denominator,
one obtains a kind of expression that we will
denote in summarized form as:
\[
\frac{
{\sf polynomial}
\big(
\big(
\Theta_{z^j\overline{z}^k\overline{w}^l}
\big)_{1\leqslant j+k+l\leqslant 6}
\big)}{
\big[\Theta_{\overline{z}}\,\Theta_{z\overline{w}}
-
\Theta_{\overline{w}}\,\Theta_{z\overline{z}}\big]^7},
\]
and hence instantly, sphericity of $M$ is characterized by:
\[
\boxed{\,
0
\equiv
{\sf polynomial}
\Big(
\big(
\Theta_{z^j\overline{z}^k\overline{w}^l}
(z,\overline{z},\overline{w})
\big)_{1\leqslant j+k+l\leqslant 6}
\Big).\,}
\]
One notices that the complex graphing function $\Theta$ is differentiated
always at leat once.

\medskip\noindent{\bf Interpretation.}
{\em Then the true thing is: after erasing the Levi determinant lying at
denominator place, if this explicit equation vanishes in some very
small neighborhood of some point: 
\[
(z_p,\overline{z}_p,\overline{w}_p) 
\,\in\, 
\square_{\rho_0}^3
\]
of the threedisc of convergence of $\Theta$,
then by the uniqueness principle for analytic functions, 
the concerned polynomial numerator:
\[
{\sf polynomial}
\big(
\big(
\Theta_{z^j\overline{z}^k\overline{w}^l}
(z,\overline{z},\overline{w})
\big)_{1\leqslant j+k+l\leqslant 6}
\big)
\]
immediately vanishes identically allover $\square_{\rho_0}^3$, so
that sphericity at one point should
freely propagate to all other Levi
nondegenerate points $q \in M \cap \square_{\rho_0}^2$.}

\medskip

Before providing rigorous details to explain the latter
assertion, a further comment is in order.

\medskip\noindent{\bf Speculative intuitive thought.}
{\em The \underline{explicit} sphericity formula brings 
the important information that denominator places are occupied
by nondegeneracy conditions, so that division is
allowed only at points where these conditions are
satisfied, but numerator places happen to be \underline{polynomial},
a computional
fact which hence enables one to jump across
degenerate points
through the `bridge-numerator' from one nondegenerate point to 
another nondegenerate point.}

\medskip

Now, the local version of Theorem~\ref{main-theorem} is as follows.
Notice that from
now on, one does not assume anymore that the origin
$0 \in M \backslash \Sigma_{\sf LD}$ be a Levi nondegenerate point.

\begin{Proposition}
With $M^3 \subset \C^2$, $\square_{\rho_0}^2$, $(z, \overline{ z}, 
\overline{w})$, 
$\varphi$, $\Theta$ as above, assuming that the real analytic subset
$\Sigma_{\sf LD}$ of Levi degenerate points is proper, 
if $M$ is spherical at one Levi nondegenerate point:
\[
p
\in
\big(M\backslash
\Sigma_{\sf LD}\big)
\cap
\square_{\rho_0}^2,
\]
then $M$ is also spherical at every other Levi nondegenerate point:
\[
q
\in
\big(M\backslash
\Sigma_{\sf LD}\big)
\cap
\square_{\rho_0}^2.
\]
\end{Proposition}

\proof
Take a (possibly much) smaller bidisc:
\[
p+\square_{\rho'}^2
\,\subset\subset\,
\square_{\rho_0}^2,
\]
with $0 < \rho' \ll \rho_0$ to be chosen below, 
and center new coordinates at:
\[
p
=
(z_p,w_p),
\]
that is to say, introduce the new translated coordinates:
\[
z'
:=
z-z_p,
\ \ \ \ \ \ \ \ \ \ \ \ \ \ \ \ \
w'
:=
w-w_p.
\]
The graphed complex equation:
\[
w
=
\Theta\big(z,\overline{z},\overline{w}\big)
\]
then becomes:
\[
w'+w_p
=
\Theta\big(z'+z_p,\,\overline{z}'+\overline{z}_p,\,\overline{w}'
+\overline{w}_p'\big).
\]
Of course, the fact that $p \in M$ reads:
\[
w_p
=
\Theta\big(z_p,\overline{z}_p,\overline{w}_p\big),
\]
and hence, in the new coordinates $(z', w')$ centered
at $p$, the equation of $M$ becomes:
\[
\boxed{\,
\aligned
w'
&
=
\Theta\big(z'+z_p,\,\overline{z}'+\overline{z}_p,\,
\overline{w}'+\overline{w}_p\big)
-
\Theta\big(z_p,\overline{z}_p,\overline{w}_p\big)
\\
&
=:
\Theta'\big(z',\overline{z}',\overline{w}'\big),
\endaligned
\,}
\]
in terms of a new graphing function $\Theta'$ that visibly satisfies:
\[
\Theta'\big(0',0',0')
=
0.
\]

\begin{Observation}
For any integers:
\[
(j,k,l)
\in
\N^3,
\]
with:
\[
j+k+l
\geqslant
1,
\]
one has:
\[
\boxed{\,\,
\Theta_{{z'}^j{\overline{z}'}^k{\overline{w}'}^l}'
(0',0',0')
=
\Theta_{z^j\overline{z}^k\overline{w}^l}
\big(z_p,\overline{z}_p,\overline{w}_p\big).\,\,
}
\]
\end{Observation}

\proof
Indeed, the constant $-\, \Theta \big( z_p, \overline{ z}_p, 
\overline{ w}_p \big)$ disappears after just a single differentiation.
\endproof

Now, the Levi nondegeneracy of $M$ at $p$ which reads
according to what precedes as:
\[
0
\neq
\big[
\Theta_{\overline{z}}\,\Theta_{z\overline{w}}
-
\Theta_{\overline{w}}\,\Theta_{z\overline{z}}
\big]
\big(z_p,\overline{z}_p,\overline{w}_p\big),
\]
reads in the new coordinates as:
\[
0
\neq
\big[
\Theta_{\overline{z}'}'\,\Theta_{z'\overline{w}'}'
-
\Theta_{\overline{w}'}'\,\Theta_{z'\overline{z}'}
\big]
(0',0',0'),
\]
which means as we know Levi nondegeneracy at $(z', \overline{ z}', 
\overline{ w}') = (0', 0', 0')$.

Precisely because in~\cite{ Merker-2010} one needs only 
this condition to hold in order to associate as was
explained above a second-order complex ordinary differential equation:
\[
w_{z'z'}'
=
\Phi'\big(
z',w'(z'),w_{z'}(z')\big)
\]
of course for some possibly very small:
\[
\vert z'\vert
<
\rho',
\ \ \ \ \ \ \ \ \ \ \ \ \ \ \ \ \ \ \ \
\vert w'\vert
<
\rho',
\]
---\,\,this is where one has to choose $\rho'$ with $0 < \rho' \ll
\rho_0$, the possible presence of rather close Levi degenerate points
being a constraint in the needed application(s) of the implicit
function theorem\,\,---, one has the impression that one can in
principle only determine whether $M$ is spherical
restrictively in such a very narrow
neighborhood $\square_{ \rho'}^2$ of $p$ in $\C^2$,
when one applies the main result of~\cite{ Merker-2010}. 

But looking just at the numerator of the
equation which expresses that $M$ is spherical
at $p$ in the coordinates $(z', w')$:
\[
0
\equiv
\frac{
\overbrace{{\sf polynomial}}^{\text{\rm same universal}
\atop
\text{\rm expression}}
\big(
\big(
\Theta_{{z'}^j{\overline{z}'}^k{\overline{w}'}^l}'
\big)_{1\leqslant j+k+l\leqslant 6}
\big)}{
\underbrace{\big[\Theta_{\overline{z}'}'\,\Theta_{z'\overline{w}'}'
-
\Theta_{\overline{w}}'\,\Theta_{z'\overline{z}'}'\big]^7}_{
\text{\rm nonvanishing at $(0',0',0')$}}},
\]
if one takes account of the above observation, one readily 
realizes that:
\[
{\sf polynomial}
\big(
\big(
\Theta_{{z'}^j{\underline{z}'}^k{\underline{w}'}^l}'
(z',\underline{z}',\underline{w}')
\big)_{1\leqslant j+k+l\leqslant 6}
\big)
=
{\sf polynomial}
\big(
\big(
\Theta_{z^j\underline{z}^k\underline{w}^l}
(z,\underline{z},\underline{w})
\big)_{1\leqslant j+k+l\leqslant 6}
\big),
\]
so that the identical vanishing of the left-hand side for:
\[
\vert z'\vert
<
\rho'
\ll
\rho_0,
\ \ \ \ \ \ \ \ \ \
\vert\underline{z}'\vert
<
\rho'
\ll
\rho_0,
\ \ \ \ \ \ \ \ \ \
\vert\underline{w}'\vert
<
\rho'
\ll
\rho_0,
\]
means the identical vanishing of the right-hand side for:
\[
\vert z-z_p\vert
<
\rho'
\ll
\rho_0,
\ \ \ \ \ \ \ \ \ \
\vert\underline{z}-\overline{z}_p\vert
<
\rho'
\ll
\rho_0,
\ \ \ \ \ \ \ \ \ \
\vert\underline{w}-\overline{w}_p\vert
<
\rho'
\ll
\rho_0,
\]
which lastly yields {\em thanks to the uniqueness principle
enjoyed by analytic functions} the
identical vanishing of the original numerator
{\em in the whole initial domain of convergence}:
\[
\boxed{\,\,
0
\equiv
{\sf polynomial}
\big(
\big(
\Theta_{z^j\underline{z}^k\underline{w}^l}
(z,\underline{z},\underline{w})
\big)_{1\leqslant j+k+l\leqslant 6}
\big)
\ \ \ \ \ \ \ \ \ \ \ \ \ 
{\scriptstyle{(\vert z\vert\,<\,\rho_0,\,\,
\vert\underline{z}\vert\,<\,\rho_0,\,\,
\vert\underline{w}\vert\,<\,\rho_0)}}.\,\,}
\]

Take now any other Levi nondegenerate point:
\[
q
\in
M\cap\square_{\rho_0}^2.
\]
The goal is to prove that $M$ is also spherical at $q$.
Center similarly new coordinates at $q = (z_q, w_q)$:
\[
z''
:=
z-z_q,
\ \ \ \ \ \ \ \ \ \ \ \ \ 
w''
:=
w-w_q.
\]
Introduce the new graphed equations:
\[
\aligned
w''
&
=
\Theta\big(z''+z_q,\,\overline{z}''+\overline{z}_q,\,\overline{w}''
+\overline{w}_q\big)
-
\Theta\big(z_q,\overline{z_q},\overline{w}_q\big)
\\
&
=:
\Theta''\big(z'',\overline{z}'',\overline{w}''\big).
\endaligned
\]
At such a point, since the Levi determinant is nonvanishing, 
one can for completeness construct the associated second-order complex
ordinary differential equations:
\[
w_{z''z''}''(z'')
=
\Phi''\big(z'',w''(z''),\,w_{z''}''(z'')\big),
\]
or question directly whether local sphericity holds near $q$
by plainly applying the main theorem of~\cite{ Merker-2010},
namely question whether the following equation holds:
\[
0\,\overset{?}{\,\equiv\,}
\frac{
\overbrace{{\sf polynomial}}^{\text{\rm again same universal}
\atop
\text{\rm expression}}
\big(
\big(
\Theta_{{z''}^j{\overline{z}''}^k{\overline{w}''}^l}''
\big)_{1\leqslant j+k+l\leqslant 6}
\big)}{
\underbrace{\big[\Theta_{\overline{z}''}''\,\Theta_{z''\overline{w}''}''
-
\Theta_{\overline{w}}''\,\Theta_{z'\overline{z}''}''\big]^7}_{
\text{\rm nonvanishing at $(0'',0'',0'')$}}}.
\]
But then by exactly the same application of the above observation,
we know that this last numerator satisfies:
\[
{\sf polynomial}
\big(
\big(
\Theta_{{z''}^j{\underline{z}''}^k{\underline{w}''}^l}''
(z'',\underline{z}'',\underline{w}'')
\big)_{1\leqslant j+k+l\leqslant 6}
\big)
=
{\sf polynomial}
\big(
\big(
\Theta_{z^j\underline{z}^k\underline{w}^l}
(z,\underline{z},\underline{w})
\big)_{1\leqslant j+k+l\leqslant 6}
\big),
\]
when:
\[
\aligned
& \ \ \ \ \ \ \ 
\vert z''\vert
<
\rho''
\ll
\rho_0,
\ \ \ \ \ \ \ \ \ \ \ \ \ \ \ \ \ \,
\vert\underline{z}''\vert
<
\rho''
\ll
\rho_0,
\ \ \ \ \ \ \ \ \ \ \ \ \ \ \ \ \ \
\vert\underline{w}'\vert
<
\rho''
\ll
\rho_0,
\\
&
\vert z-z_q\vert
<
\rho''
\ll
\rho_0,
\ \ \ \ \ \ \ \ \ \
\vert\underline{z}-\overline{z}_q\vert
<
\rho''
\ll
\rho_0,
\ \ \ \ \ \ \ \ \ \
\vert\underline{w}-\overline{w}_q\vert
<
\rho''
\ll
\rho_0,
\endaligned
\]
and since we already know that the latter right-hand side vanishes,
according to the last boxed equation, we conclude that $M$ is indeed
spherical at $q$.
\endproof

To finish the proof of Theorem~1.1 in the case $n = 1$, it remains
only to {\em globalize} this local propagation of sphericity. One
does this by means of standard arguments which consist to pick up one
Levi nondegenerate point $p^\sim \in \square_{ \rho_0}^2$ (possibly
close to the boundary of the bidisc!), to center some affine
coordinates at $p^\sim$, to use local real analytic equations for $M$
expanded in some Taylor series which converge in some other bidisc
$\square_{ \rho_0^\sim}^2$ centered
at $p^\sim$, and to apply the same reasonings as above
to propagate sphericity from $p^\sim$ to any other Levi nondegenerate
point $q^\sim \in M \cap \square_{\rho_0^\sim}^2$. By connectedness of
$M$, one concludes.

\section{Proof in $\C^{n+1}$ ($n \geqslant 2$)} 
\label{C-n-1}
\HEAD{\ref{C-n-1}.~Proof in $\C^{n+1}$ ($n \geqslant 2$)}{
Jo\"el Merker, D\'epartement de Math\'ematiques d'Orsay, Paris, France}

We briefly summarize the quite similar arguments,
relying upon~\cite{ Merker-2009}.

The Levi determinant becomes:
\[
\Delta
:=
\left\vert
\begin{array}{cccc}
\Theta_{\overline{z}_1} & \cdots & \Theta_{\overline{z}_n} 
& \Theta_{\overline{w}}
\\
\Theta_{z_1\overline{z}_1} & \cdots & \Theta_{z_1\overline{z}_n} 
& \Theta_{z_1\overline{w}}
\\
\cdot\cdot & \cdots & \cdot\cdot & \cdot\cdot
\\
\Theta_{z_n\overline{z}_1} & \cdots & \Theta_{z_n\overline{z}_n} 
& \Theta_{z_n\overline{w}}
\end{array}
\right\vert.
\]
It is nonzero at one point:
\[
p
=
\big(z_{1p},\dots,z_{np},w_p\big)
\in
M
\]
if and only if $M$ is Levi nondegenerate at $p$, and also, if and only
if one can associate to $M$ a completely integrable system of
second-order partial differential equations:
\[
w_{z_{k_1}z_{k_2}}(z)
=
\Phi_{k_1,k_2}
\big(
z,\,w(z),\,w_{z_1}(z),\dots,w_{z_n}(z)
\big)
\ \ \ \ \ \ \ \ \ \ \ \ \ \ \ \ \ \
{\scriptstyle{(1\,\leqslant\,k_1,\,\,k_2\,\leqslant\,n)}}.
\]
Hachtroudi (\cite{ Hachtroudi-1937}) established that such a system is
pointwise equivalent to:
\[
w_{z_{k_1}'z_{k_2}'}'(z')
=
0
\ \ \ \ \ \ \ \ \ \ \ \ \ \ \ \ \ \
{\scriptstyle{(1\,\leqslant\,k_1,\,\,k_2\,\leqslant\,n)}}
\]
if and only if:
\[
\footnotesize
\aligned
0
&
\equiv
\frac{\partial^2\Phi_{k_1,k_2}}{\partial w_{z_{\ell_1}}w_{z_{\ell_2}}}
-
\\
&
\ \ \ \ \
-
{\textstyle{\frac{1}{n+2}}}\,
\sum_{\ell_3=1}^n
\left(
\delta_{k_1,\ell_1}
\frac{\partial^2\Phi_{\ell_3,k_2}}{
\partial w_{z_{\ell_3}}\partial w_{z_{\ell_2}}}
+
\delta_{k_1,\ell_2}
\frac{\partial^2\Phi_{\ell_3,k_2}}{
\partial w_{z_{\ell_1}}\partial w_{z_{\ell_3}}}
+
\delta_{k_2,\ell_1}
\frac{\partial^2\Phi_{k_1,\ell_3}}{
\partial w_{z_{\ell_3}}\partial w_{z_{\ell_2}}}
+
\delta_{k_2,\ell_2}
\frac{\partial^2\Phi_{k_1,\ell_3}}{
\partial w_{z_{\ell_1}}\partial w_{z_{\ell_3}}}
\right)
+
\\
&
\ \ \ \ \
+
{\textstyle{\frac{1}{
(n+1)(n+2)}}}
\big[
\delta_{k_1,\ell_1}\delta_{k_2,\ell_2}
+
\delta_{k_2,\ell_1}\delta_{k_1,\ell_2}
\big]
\sum_{\ell_3=1}^n\,\sum_{\ell_4=1}^n\,
\frac{\partial^2\Phi_{\ell_3,\ell_4}}{
\partial w_{z_{\ell_3}}\partial w_{z_{\ell_4}}}
\ \ \ \ \ \ \ \ \ \ \ \ \ \ \ \ \ \ \ \ \ \ \ \
\begin{array}{c}
{\scriptstyle{(1\,\leqslant\,k_1,\,\,k_2\,\leqslant\,n)}}
\\
{\scriptstyle{(1\,\leqslant\,\ell_1,\,\,\ell_2\,\leqslant\,n)}}
\end{array}.
\endaligned
\]

When one does apply Hachtroudi's results to CR geometry (instead of
Chern-Moser's, which is up to now not sufficiently explicit to be
applied), the signature of the Levi forms disappears for the following
reason.

The infinite-dimensional local Lie (pseudo-)group of biholomorphic
transformations:
\[
\aligned
\big(z_1,\dots,z_n,w\big)
&
\,\longmapsto\,
\big(z_1',\dots,z_n',w'\big)
\\
&\ \ \ \ \,
=
\big(z_1'(z_\bullet,w),\dots,z_n'(z_\bullet,w),w'(z_\bullet,w)\big)
\endaligned
\]
acts simultaneously on $(z_\bullet, w)$-variables and on $(\overline{
z}_\bullet, \overline{ w})$-variables as:
\[
\aligned
\big(\overline{z}_1,\dots,\overline{z}_n,\overline{w}\big)
&
\,\longmapsto\,
\big(
\overline{z}_1',\dots,\overline{z}_n',\overline{w}'\big)
\\
&\ \ \ \ \,
=
\big(\overline{z}_1'(
\overline{z}_\bullet,
\overline{w}),\dots,
\overline{z}_n'(
\overline{z}_\bullet,
\overline{w}),
\overline{w}'(
\overline{z}_\bullet,
\overline{w})\big).
\endaligned
\]

But when one passes to the extrinsic complexification,
one replaces $(\overline{ z}, \dots, \overline{ z}_n, 
\overline{ w})$-variables by new independent variables:
\[
\big(\underline{z}_1,\dots,\underline{z}_n,\underline{w}\big),
\]
considered as the {\em constants} of integration for the system of
partial differential equations. Hence the local Lie (pseudo)-group
considered by Hachtroudi
becomes enlarged as the group of transformations:
\[
\big(z_\bullet,w,\,\underline{z}_\bullet,\underline{w}\big)
\,\longmapsto\,
\Big(
\text{\footnotesize\sf holomorphic map}
\big(z_\bullet,w\big),\,\,
\text{\footnotesize\sf other holomorphic map}
\big(\underline{z}_\bullet,\underline{w}\big)
\Big)
\]
in which the transformations on the `constant-of-integration'
variables $(\underline{ z}_1, \dots, \underline{ z}_n, \underline{
w})$ becomes completely dis-coupled from the group of
transformations on the true variables $(z_1, \dots, z_n, w)$. By
definition ({\em cf.}~the explanations 
in~\cite{ Merker-2010}), transformations on
differential equations, when viewed in the space of solutions, are
always of this general form.

It is then clear that all {\em complexified} Heisenberg $(k, n-k)$
pseudospheres:
\[
w
=
\underline{w}
+
2i\big(
-\,z_1\underline{z}_1
-\cdots-
z_k\underline{z}_k
+
z_{k+1}\underline{z}_{k+1}
+\cdots+
z_n\underline{z}_n
\big),
\] 
become all pairwise equivalent through such transformations, because one is
allowed to replace $\underline{ z}_1, \dots, \underline{ z}_k$ by $-\,
\underline{ z}_1, \dots, - \underline{ z}_k$ without touching $z_1,
\dots, z_k$; even the factor $i$ can be erased:
\[
w
=
\underline{w}
+
z_1\underline{z}_1
+\cdots+
z_k\underline{z}_k
+
z_{k+1}\underline{z}_{k+1}
+\cdots+
z_n\underline{z}_n
\] 
Therefore, when passing to systems of partial differential
equations associated to CR manifolds, {\em Levi form
signatures drop}.

Consequently, when one applies the main theorem of~\cite{
Merker-2009}, according to which a Levi nondegenerate
$M \subset \C^{ n+1}$ {\em having given Levi form signature
$(k, n-k)$} is pseudo-spherical if and only if
(notation same as in~\cite{ Merker-2009}) its
Hachtroudi system is equivalent to:
\[
w_{z_{k_1}'z_{k_2}'}'(z')
=
0
\ \ \ \ \ \ \ \ \ \ \ \ \ \ \ \ \ \
{\scriptstyle{(1\,\leqslant\,k_1,\,\,k_2\,\leqslant\,n)}},
\]
and moreover, if and only if\,\,---\,\,after translating back to 
the graphing function $\Theta$ the explicit 
condition of Hachtroudi\,\,---\,\,the following
identical vanishing property holds:
\[
\scriptsize
\aligned
0
&
\equiv
\frac{1}{\Delta^3}
\bigg[\,
\sum_{\mu=1}^{n+1}\,\sum_{\nu=1}^{n+1}
\bigg[
\Delta_{[0_{1+\ell_1}]}^\mu
\cdot
\Delta_{[0_{1+\ell_2}]}^\nu
\bigg\{
\Delta
\cdot
\frac{\partial^4\Theta}{
\partial z_{k_1}\partial z_{k_2}
\partial\overline{t}_\mu\partial\overline{t}_\nu}
-
\sum_{\tau=1}^{n+1}\,
\Delta_{[\overline{t}^\mu\overline{t}^\nu]}^\tau
\cdot
\frac{\partial^3\Theta}{
\partial z_{k_1}\partial z_{k_2}\partial\overline{t}^\tau}
\bigg\}
-
\\
&
-
{\textstyle{\frac{\delta_{k_1,\ell_1}}{n+2}}}\,
\sum_{\ell_3=1}^n\,
\Delta_{[0_{1+\ell_3}]}^\mu
\cdot
\Delta_{[0_{1+\ell_2}]}^\nu
\bigg\{
\Delta
\cdot
\frac{\partial^4\Theta}{
\partial z_{\ell_3}\partial z_{k_2}
\partial\overline{t}_\mu\partial\overline{t}_\nu}
-
\sum_{\tau=1}^{n+1}\,
\Delta_{[\overline{t}^\mu\overline{t}^\nu]}^\tau
\cdot
\frac{\partial^3\Theta}{
\partial z_{\ell_3}\partial z_{k_2}\partial\overline{t}^\tau}
\bigg\}
-
\\
&
-
{\textstyle{\frac{\delta_{k_1,\ell_2}}{n+2}}}\,
\sum_{\ell_3=1}^n\,
\Delta_{[0_{1+\ell_1}]}^\mu
\cdot
\Delta_{[0_{1+\ell_3}]}^\nu
\bigg\{
\Delta
\cdot
\frac{\partial^4\Theta}{
\partial z_{\ell_3}\partial z_{k_2}
\partial\overline{t}_\mu\partial\overline{t}_\nu}
-
\sum_{\tau=1}^{n+1}\,
\Delta_{[\overline{t}^\mu\overline{t}^\nu]}^\tau
\cdot
\frac{\partial^3\Theta}{
\partial z_{\ell_3}\partial z_{k_2}\partial\overline{t}^\tau}
\bigg\}
-
\\
&
-
{\textstyle{\frac{\delta_{k_2,\ell_1}}{n+2}}}\,
\sum_{\ell_3=1}^n\,
\Delta_{[0_{1+\ell_3}]}^\mu
\cdot
\Delta_{[0_{1+\ell_2}]}^\nu
\bigg\{
\Delta
\cdot
\frac{\partial^4\Theta}{
\partial z_{k_1}\partial z_{\ell_3}
\partial\overline{t}_\mu\partial\overline{t}_\nu}
-
\sum_{\tau=1}^{n+1}\,
\Delta_{[\overline{t}^\mu\overline{t}^\nu]}^\tau
\cdot
\frac{\partial^3\Theta}{
\partial z_{k_1}\partial z_{\ell_3}\partial\overline{t}^\tau}
\bigg\}
-
\\
&
-
{\textstyle{\frac{\delta_{k_2,\ell_2}}{n+2}}}\,
\sum_{\ell_3=1}^n\,
\Delta_{[0_{1+\ell_1}]}^\mu
\cdot
\Delta_{[0_{1+\ell_3}]}^\nu
\bigg\{
\Delta
\cdot
\frac{\partial^4\Theta}{
\partial z_{k_1}\partial z_{\ell_3}
\partial\overline{t}_\mu\partial\overline{t}_\nu}
-
\sum_{\tau=1}^{n+1}\,
\Delta_{[\overline{t}^\mu\overline{t}^\nu]}^\tau
\cdot
\frac{\partial^3\Theta}{
\partial z_{k_1}\partial z_{\ell_3}\partial\overline{t}^\tau}
\bigg\}
+
\\
&
\ \ \ \ \
+
{\textstyle{\frac{1}{(n+1)(n+2)}}}
\cdot
\big[
\delta_{k_1,\ell_1}\delta_{k_2,\ell_2}
+
\delta_{k_2,\ell_1}\delta_{k_1,\ell_2}
\big]
\cdot
\\
&
\ \ \ \ \
\cdot
\sum_{\ell_3=1}^n\,\sum_{\ell_4=1}^n\,
\Delta_{[0_{1+\ell_3}]}^\mu
\cdot
\Delta_{[0_{1+\ell_4}]}^\nu
\bigg\{
\Delta
\cdot
\frac{\partial^4\Theta}{
\partial z_{\ell_3}\partial z_{\ell_4}
\partial\overline{t}_\mu\partial\overline{t}_\nu}
-
\sum_{\tau=1}^{n+1}\,
\Delta_{[\overline{t}^\mu\overline{t}^\nu]}^\tau
\cdot
\frac{\partial^3\Theta}{
\partial z_{\ell_3}\partial z_{\ell_4}\partial\overline{t}^\tau}
\bigg\}
\bigg],
\\
&
\ \ \ \ \ \ \ \ \ \ \ \ \ \ \ \ \ \ \ \ \ \ \ \ \ \ \ \ \ \ 
\ \ \ \ \ \ \ \ \ \ \ \ \ 
(1\,\leqslant k_1,\,k_2\,\leqslant\,n;\,\,\,
(1\,\leqslant \ell_1,\,\ell_2\,\leqslant\,n),
\endaligned
\] 
one can reason exactly as in the preceding section 
for $M^3 \subset \C^2$\,\,---\,\,noticing 
that the denominator is similarly $\frac{ 1}{ \Delta^3}$, 
noticing that the numerator is similarly polynomial
in the partial derivatives of the graphing function
$\Theta$\,\,---, but
when one jumps from a Levi nondegenerate
point $p \in M \cap \square_{\rho_0}^{n+1}$ to
another Levi nondegenerate point $q \in M \cap 
\square_{ \rho_0}^{n+1}$, from the property 
of local equivalence at $q$ to:
\[
w_{z_{k_1}''z_{k_2}''}''(z'')
=
0
\ \ \ \ \ \ \ \ \ \ \ \ \ \ \ \ \ \
{\scriptstyle{(1\,\leqslant\,k_1,\,\,k_2\,\leqslant\,n)}},
\]
one can only conclude that the complexification
of $M$ near $q$ is equivalent near $q$ to:
\[
w''
=
\underline{w}''
+
z_1''\underline{z}_1''
+\cdots+
z_k''\underline{z}_k''
+
z_{k+1}''\underline{z}_{k+1}''
+\cdots+
z_n''\underline{z}_n''
\] 
so that the Levi form signature can in 
principle change\,\,---\,\,and really does in Example~6.3
of~\cite{ Kossovskiy-Shafikov-2012}\,\,---\,\,through
Levi degenerate points.
\qed

\vfill\end{document}